\documentclass[preprint,12pt]{elsarticle}

\usepackage{graphicx,color}

\usepackage{amssymb}
\usepackage{amsthm}

\usepackage{amsmath,mathrsfs,bm,subfigure,verbatim}

\journal{Journal of Computational and Applied Mathematics}

\newtheorem{thm}{Theorem}

\newdefinition{defn}{Definition}
\newdefinition{rmk}{Remark}
\newdefinition{alg}{Algorithm}
\newdefinition{exmp}{Example}
\newproof{pf}{Proof}

\begin{document}

\begin{frontmatter}




\title{Fitzpatrick Algorithm  for Multivariate Rational Interpolation}
\tnotetext[fund]{This work was supported in part by the National
Grand Fundamental Research 973 Program of China (No. 2004CB318000).}


\author{Peng Xia }
\author{Shugong Zhang\corref{cor1}}
\cortext[cor1]{Corresponding author} \ead{sgzh@mail.jlu.edu.cn}
\author{Na Lei}
\address{School of Mathematics,
Key Lab. of Symbolic Computation and Knowledge Engineering
\textup{(}Ministry of Education\textup{)}, Jilin University,
Changchun 130012, PR China}

\begin{abstract}
In this paper, we first apply the Fitzpatrick algorithm to
osculatory rational interpolation. Then based on Fitzpatrick
algorithm, we present a Neville-like algorithm for Cauchy
interpolation. With this algorithm, we can determine the value of
the interpolating function at a single point without computing the
rational interpolating function.

\end{abstract}


\end{frontmatter}


\section{Introduction}
\label{sec:int} 
Interpolation is an important method in numerical approximation.
Rational functions sometimes are superior to polynomial for
interpolating data because they can achieve more accurate
approximations with the same amount of computation \cite{Key17}. In
addition, rational interpolants have a natural way of interpolating
poles whereas polynomial interpolants do not. So how to solve the
problem of rational interpolation is what people have been
concerning.

There are rich literatures on the univariate Cauchy interpolation
and osculatory rational interpolation problem, such as \cite{Key14,
Key15, Key12, Key38, Key11, Key9, Key10,
 Key13}.
For multivariate rational interpolation \cite{Key14, Key15, Key12,
Key2, Key16} gave some results about bivariate cases and the authors
assume the interpolation nodes are on rectangular girds.
\cite{Key36} computed rational interpolation over pyramid-typed
grids in $\mathbb{R}^3$ by branched continued fractions. When the
interpolation data are scattered, \cite{Key6} presented a fast
solver for the linear block Cauchy-Vandermonde system that
translates the interpolation conditions. \cite{Key22} used the
theory of algebraic geometry to study the minimal multivariate
rational interpolation.

One of the main problems of rational interpolation is the
parametrization of all solutions of a given degree of complexity.
Based on Euclidean algorithm, \cite{Key35} investigated a general
frame work which lead to a parametrization of all rational
interpolation functions. \cite{Key31} considered the set $M=\{(a,b):
a\equiv bh \mod x^{2t}\}$ of all solutions of the key equation for
alternant codes, and give the Fitzpatrick algorithm. \cite{Key32}
extended the Fitzpatrick algorithm to determining a parametrization
of all minimal complexity rational functions $a(x)/b(x)$
interpolating an arbitrary sequence of points, and complexity is
measured in terms of $\max\{\deg ( a(x))$, $\deg(b(x))+\xi\}$, where
$\xi$ is a given integer. \cite{Key34} presented an algorithm to
seeking the Gr\"{o}bner basis for the  solution of polynomial
congruences in one or more variables. \cite{Key33} generalized the
work in \cite{Key34}, and got a general algorithm applicable to a
wide range of constrained interpolation.

In this paper we apply the Fitzpatrick algorithm which appears in
\cite{Key33} to multivariate osculatory rational interpolation, and
get the parametric solution of the multivariate osculatory rational
interpolation function $ r(x)= a(X)/b(X)$. Based on Fitzpatrick
algorithm, we present a Neville-like algorithm for multivariate
Cauchy interpolation.

The rest of the paper is organized as follows. In Section 2, we
present the Fitzpatrick algorithm. In Section 3, we apply the
Fitzpatrick algorithm to seek the weak solution $( a(X),b(X))$ of
multivariate osculatory rational interpolation. In Section 4, based
on Fitzpatrick algorithm, we give a Neville-like algorithm for
multivariate Cauchy interpolation. With this algorithm, we can
determine the value of the interpolating function at a single point
without computing the rational interpolating function.

\section{Fitzpatrick algorithm}\label{FA}

Fitzpatrick algorithm, also called FGLM-like algotithm, is
applicable to coding theory, Pad$\acute{e}$ approximation, partial
realization, interpolation, and other fields. \cite{Key30} described
these literatures from a historical point of view, started from
\cite{Key31}, and covered recent developments for list decoding. For
further details, please refer to \cite{Key30} or \cite{Key1}. Now we
introduce the structure of Fitzpatrick algorithm.

Let $\mathbb{F}$ be a field, $\mathcal
{P}=\mathbb{F}[x_1,\ldots,x_n]$ be a polynomial ring and $d\geq 1$
be a natural number.

We denote by $\{\vec {e}_1,\ldots,\vec {e}_d\}$ the canonical basis
of $\mathcal{P}^d$. Any term in $\mathcal{P}^d$ is of the form
$m=\phi \vec {e}_k$, $1\leq k \leq d$, where $\phi$ is a term in
$\mathcal{P}$, and the set of terms in $\mathcal{P}^d$ is denoted by
$\mathcal {T}^{(d)}$. Let $\prec$ be a term order. For each $ \vec
f=\underset{\tau \in \mathcal {T}^{(d)}}{\sum} c(\vec f,\tau)\tau \
\in \mathcal{P}^d$, its support is
$$supp(\vec f):=\{\tau \in \mathcal{T}^d : c(\vec f,\tau)\neq 0\},$$
its $leading$ $term$ is $\textbf{LT}
(\vec
f):=\max_{\prec}(supp(\vec f))$, its $leading$ $coefficient$ is
$\textbf{LC}
(\vec f):=c(\vec f,\textbf{LT}
(\vec f))$
and $leading$ $monomial$ is $\textbf{LM}
(\vec
f):=\textbf{LC}
(\vec f)\textbf{LT}
(\vec f)$.


Let $M_k$, $M_{k+1}$ be submodules of a $\mathcal{P}$-module $M$,
with $M_k \supseteq M_{k+1}$, such that, for each $s$, $1\leq s\leq
n$, there exists $\beta_s \in \mathbb{F}$ satisfying
\begin{equation}\label{e2.1}
    (x_s-\beta_s)M_{k}\subseteq M_{k+1}.
\end{equation}
For each $k$, there exists an $\mathbb{F}$-homomorphism
\begin{equation}\label{e2.2}
    \theta_k : M_k\longrightarrow \mathbb{F}
\end{equation}
with $\ker(\theta_k)=M_{k+1}$.

In \cite{Key33}, they described Fitzpatrick algorithm in the
following theorem.
\begin{thm}\textup{\cite{Key33}}\label{th1}
Let $M$ be a $\mathcal{P}$-module and let $M_{k}\supseteq M_{k+1}$
be submodules of $M$ satisfying $(\ref{e2.1})$ and $(\ref{e2.2})$
for suitable $\beta_s$, $\theta_k$. Let $H
:\mathcal{P}^d\longrightarrow M$ be an $\mathbb{F}$-linear function
such that for each $s$, $1\leq s\leq n$, there exists $\gamma_s\in
\mathbb{F}$ satisfying
$$H(x_s \vec {b})=(x_s+\gamma_s)H(\vec {b})$$
for all $\vec {b}=(b_1,\ldots,b_d)\in \mathcal{P}^d$. Let
$S\subseteq \mathcal{P}^d$ be a submodule satisfying
\begin{equation}\label{e2.3}
    H(\vec {b})\equiv 0 \mod M_{k},\  \forall \vec {b}\in S
\end{equation}
and let $S'\subseteq S$ be the set of elements satisfying
\begin{equation}\label{e2.4}
    H(\vec {b})\equiv 0 \mod M_{k+1}
\end{equation}
Then $S'$ is a submodule of $\mathcal{P}^d$.

If we have obtained an ordered minimal Gr\"{o}bner basis
$G=\{G[1],\ldots,$ $G[|G|]\}$ of $S$ with respect to a term order
$\prec$, then a Gr\"{o}bner basis $G'$ of $S'$ with respect to
$\prec$ can be constructed as follows:

Define $\alpha_j=\theta_k(H(G[j]))$, for $1\leq j\leq |G|$.

If $\alpha_j=0$ for all $j$ then

$\qquad$ $G'=G$

otherwise

$\qquad$ $j^{*}$= the least $j$ for which $\alpha_j\neq 0$

$\qquad$ $G_1:=\{G[j]:j<j^{*}\}$

$\qquad$ $G_2:=\{(x_s-(\beta_s+\gamma_s))G[j^{*}]: 1\leq s\leq n\}$

$\qquad$ $G_3:=\{G[j]-(\alpha_j/\alpha_{j^{*}}) G[j^{*}]: j>j^{*}\}$

$\qquad$ $G':=G_1\cup G_2 \cup G_3$.

\end{thm}

In the following section we will use the Fitzpatrick algorithm to
compute osculatory rational interpolation functions.

\section{Osculatory rational interpolation and Fitzpatrick algorithm}\label{ORIFA}

Now, we introduce some notations. Let $\alpha$
$=(\alpha_1,\cdots,\alpha_n)\in \mathbb{N}^n$. We define a
differential operator $D^\alpha$ by
$$D^{\alpha} = \cfrac{1}{\alpha_1!\cdots \alpha_n!}\ \cfrac{\partial^{\alpha_1+\cdots+\alpha_n}}{\partial x_1^{\alpha_1}\cdots\partial x_n^{\alpha_n}}
\overset{\Delta}{=}\cfrac{1}{\alpha !}\ \cfrac{\partial^{\alpha_1+\cdots+\alpha_n}}{\partial x_1^{\alpha_1}\cdots\partial x_n^{\alpha_n}}.$$
A subset $\mathcal {A}\subset \mathbb{N}^n$ is  called a $lower$
$set$, if it is closed under the division order, that is, if
$\alpha=(\alpha_1,\ldots,\alpha_n) \in \mathcal {A}$ then $\beta \in
\mathcal {A}$ for all $ 
\beta=(\beta_1,\ldots,\beta_n)$ with $\beta_i\leq \alpha_i, i=1,\ldots, n$. 

The multivariate osculatory rational interpolation problem can be
stated as follows:

Given a set of $L$ distinct points $\{Y_1,\ldots,Y_L \}$ in space
$\mathbb{F}^n$.  Point $Y_i$ has multiplicity defined by the lower
set $\mathcal {A}_i$, and the corresponding values $\{
f_i^{(\alpha)}\in \mathbb{F} : \forall \alpha \in \mathcal
{A}_i,i=1,\ldots,L\}$. Construct a
  rational interpolation function
\begin{equation*}
  r(X)=\cfrac{  a(X)}{b(X)},
\end{equation*}
such that
\begin{equation*}
D^\alpha   r(X)\big|_{Y_i}=  f_i^{(\alpha)}~,~\forall \alpha \in
\mathcal {A}_i,~i=1,\ldots,L,
\end{equation*}
where $  a\in \mathcal {P}$, $b\in \mathcal {P}$, $  a\neq 0$,
$b(Y_i)\neq 0$ for all $i$.

 We call this problem  multivariate osculatory
rational interpolation, and $  r(X)$ multivariate osculatory
rational interpolation function.

\subsection{Weak interpolation}

From the definition of  multivariate osculatory rational
interpolation we know that the equivalent definition is: the Taylor
series expansion of $ r(X)$ at the point $X=Y_i$ satisfies
\begin{eqnarray*}
    r(X) &=& \sum_{\alpha_{i,j}\in\mathcal
{A}_i}(X-Y_i)^{\alpha_{i,j}} f_i^{(\alpha_{i,j})}+ \cdots
.
\end{eqnarray*}


Let $s_i=\sharp \mathcal {A}_i$, $i=1,\ldots,L$,
 $N=\overset{L}{\underset{i=1}\sum} s_i$. For each point $Y_i$ and
the lower set $\mathcal {A}_i$, define polynomial  $  h_i$
\begin{equation*}
  h_i := \underset{\alpha_{i,j} \in \mathcal
{A}_i}{\sum}(X-Y_i)^{\alpha_{i,j}} f_i^{(\alpha_{i,j})}.
\end{equation*}

For each $\mathcal {A}_i$, rearrange the elements of $\mathcal
{A}_i$, such that each subset $\mathcal
{A}_{i,j}=\{\alpha_{i,0},\ldots$, $\alpha_{i,j}\}, 0\leq j\leq
s_i-1$, is still a lower set. In particular, $\mathcal
{A}_{i,s_i-1}=\mathcal {A}_{i}$.

Denote the ideal $I((Y_i,\mathcal{A}_{i,j}))=\{p\in \mathcal {P} :
D^\alpha p(X)\big|_{Y_i}=0,\forall \alpha \in \mathcal{A}_{i,j}\}$
by $I_{i,j}$ and call $I_{i,j}$ the $vanishing$ $ideal$ of
$(Y_i,\mathcal{A}_{i,j})$.

\begin{defn}$(\textup{Weak}$ $\textup{interpolation})$ A pair $(  a,b)\in \mathcal{P}^{2}$ is called
a $weak$ $interpolation$ for multivariate osculatory rational
interpolation problem if
\begin{equation*}
  a\equiv b   h_i \mod I_{i,s_i-1},i=1,\ldots,L.
\end{equation*}
\end{defn}
Define $(  a,b)+(  c,d)=(  a+  c,b+d)$, $d(  a,b)=(d  a,db)$. Thus
$M=\{(  a,b):  a\equiv b  h_i \mod I_{i,s_i-1}, i=1,\ldots,L\}$ is a
$\mathcal{P}$-submodule.

If $\{(  a_1,b_1),\ldots,(  a_t,b_t)\}$ is a Gr\"{o}bner basis of
the submodule $M$, then any pair $(a,b)$ with the form
$$(  a,b)=c_1(  a_1,b_1)+\cdots+c_t(  a_t,b_t)$$
is a weak interpolation, where $c_j\in \mathcal {P}$ \
$(j=1,\ldots,t)$ are free parameters. Choose $c_j$ properly such
that $b(Y_i)\neq 0,\ i=1,\ldots,L$, then we can get the
interpolation function
$$\cfrac{  a(X)}{b(X)}=\cfrac{c_1  a_1+\cdots +c_t  a_t}{c_1b_1+\cdots+c_tb_t}.$$

\subsection{Fitzpatrick algorithm\_RI}

The aim of this subsection is to apply the Fitzpatrick algorithm to
compute osculatory rational interpolation.

\begin{defn}$(\textup{order} \prec_{\xi} )$
\begin{enumerate}
  \item We say $X^{\alpha}(  1,0) \prec_{\xi} X^{\beta}(  1,0) $
  if $|\alpha|<|\beta|$, or
  $|\alpha|=|\beta|$\  and $X^{\alpha}\prec_{lex} X^{\beta},$
  \item We say $X^{\alpha}(  1,0) \prec_{\xi} X^{\beta}(
  0,1)$ if $|\alpha|\leq|\beta|+\xi,$
  \item We say $X^{\alpha}(  0,1) \prec_{\xi} X^{\beta}(  0,1) $ if
  $|\alpha|<|\beta|$, or $|\alpha|=|\beta|$ and $X^{\alpha}\prec_{lex}
  X^{\beta},$
\end{enumerate}
where $\prec_{lex}$ is the lexicographic order on $\mathcal{P}$, and
$\xi$ is a given integer.
\end{defn}
It is easy to check that the order $\prec_{\xi}$ is a monomial order
on $\mathcal{P}^{2}$.

For each $\mathcal {A}_{i,j}$, $1\leq i\leq L$, $0\leq j\leq s_i-1$,
define the $congruent$ $equation$ as $$  a\equiv b  h_i \mod
I_{i,j},$$ where $I_{i, j}$ is the vanishing ideal of $(Y_i,
\mathcal{A}_{i, j})$.

Define an order $<$ on the lower sets $\{\mathcal {A}_{i,j}, 1\leq
i\leq L, 0\leq j\leq s_i-1\}$ such that $\mathcal {A}_{i_1,
j_1}<\mathcal {A}_{i_2, j_2}$ if and only if $i_1<i_2$ ,\ or
$i_1=i_2=i$\ and $\mathcal {A}_{i, j_1}\subset \mathcal {A}_{i,
j_2}$ for $j_1<j_2$.

Consequently, an order on the congruent equations is induced:

$  a\equiv b  h_i \mod I_{i_1,j_1}\ <\   a\equiv b  h_i \mod
I_{i_2,j_2}$ if and only if $\mathcal {A}_{i_1, j_1}<\mathcal
{A}_{i_2, j_2}$.

Now we can establish a one to one correspondence between index $k$
and $(i,j)$.  Define a sequence of submodules $M_k$,\
$k=0,\ldots,N$, where $M_0=\mathcal {P}^{2}$, $M_k$\ is the set of
common solutions of the first $k$ congruent equations
$$M_k=\{(  a,b)\in M_{k-1}:  a\equiv b  h_{i_k} \mod I_{i_k,j_k},
\}, k=1,\ldots,N.$$

Obviously $M=M_N=\{(  a,b):  a\equiv b  h_i \mod I_{i,s_i-1},
i=1,\ldots,L\}$, and $M_0\supseteq  M_1\supseteq  \cdots \supseteq
M_N$.

Fix an order $\prec_{\xi}$ we compute the minimal Gr\"{o}bner basis
$\mathcal {G}_{N}$ of $M_{N}$\ recursively. It is easy to see that
$\{ ( 1,0)$,  $( 0,1)\} $ is a Gr\"{o}bner basis of $M_0$.\ We
compute the Gr\"{o}bner basis $\mathcal {G}_{k+1}$ of $M_{k+1}$\
through the minimal Gr\"{o}bner basis $\mathcal {G}_k$ of $M_k$.

Let $\mathcal {G}_k=\{(  a_1,b_1),\ldots,( a_{m_{k}},b_{m_{k}})\}$
be the minimal Gr\"{o}bner basis of $M_k$, and the $(k+1)$-th
congruent equation be $  a\equiv b  h_l \mod I_{l,k_l}$. Then
\begin{enumerate}
  \item if $k_l=0$,\ that is $\mathcal {A}_{l,k_l}=\mathcal
  {A}_{l,0}=\{\alpha_{l,0}=  0\}$,\
  then $$b h_{l}-a\equiv \nu \equiv \nu (X-Y_l)^{\alpha_{l,0}} \mod
  I_{l,0}.$$
  \item if $k_l\neq 0$,\ that is $\mathcal {A}_{l,k_l}\big{\backslash} \mathcal
{A}_{l,k_l-1}=\{\alpha_{l,k_l}\}$,\ then
  $$bh_{l}-a \equiv \nu' (X-Y_l)^{\alpha_{l,k_l}}\mod I_{l,k_l}.$$
\end{enumerate}
Therefore for any $(  a,b)\in M_k$, we have\ $$bh_{l}-a \equiv \nu
(X-Y_l)^{\alpha_{l,k_l}}\mod I_{l,k_l},$$ and $(a,b)\in M_{k+1}$ if
and only if $\nu= 0$.

Define an $\mathbb{F}$-homonorphism
\begin{eqnarray*}
  ~~\theta &:& M_k~~ \longrightarrow \mathbb{F} \\
    & & (a,b) \longmapsto \nu
\end{eqnarray*}
Obviously, $\ker(\theta)=M_{k+1}$, $(x_s-y_{l,s})M_k\subseteq
M_{k+1}$. We define $H$: $H\big((a,b)\big)=(a,b)$. Then for any
$(a,b)\in \mathcal {P}^2$,
$H\big((x_s-y_{l,s})(a,b)\big)=(x_s-y_{l,s})H\big((a,b)\big)$. Let
$S=M_k$, $S'=M_{k+1}$.

With the definitions above, we can give the form of
Fitzpatrick algorithm for multivariate osculatory rational
interpolation.

\textbf{Fitzpatrick algorithm\_RI :}(Using the minimal Gr\"{o}bner
basis $\mathcal{G}_k$, compute the minimal Gr\"{o}bner basis of
$M_{k+1}$)
\\

     \textbf{Input:} the minimal Gr\"{o}bner
basis $\mathcal{G}_k=\{(  a_1,b_1),\ldots,(
     a_{m_{k}},b_{m_{k}})\};$

     \textbf{Output:} the minimal Gr\"{o}bner
basis  of $M_{k+1}$, that is $\mathcal{G}_{k+1}$;
\\

     $\textup{Rearrange}$ $\textup{the}$ $\textup{elements}$ $\textup{of}$ $\mathcal{G}_{k}$ $\textup{such}$
     $\textup{that}$ $\textup{LT}(  a_1,b_1)\prec_{\xi}\cdots\prec_{\xi}\textup{LT}(  a_{m_{k}},b_{m_{k}})$;

     $\textup{\textbf{for}}$ ${t}$  $\textup{from}$ $1$ $\textup{to}$ $m_{k}$ $\textup{\textbf{do}}$

     $$b_{t}h_{l}-a_{t}\equiv \nu_{t} (X-Y_l)^{\alpha_l,k_l} \mod
     I_{l,k_l};$$

     $\textup{\textbf{end}}$ $\textup{\textbf{do}}$;

     \textbf{if} $\nu_{t}=0$ for all $t$ \textbf{then}

     $\hspace{1cm}$ $\mathcal{G}_{k+1}:=\mathcal{G}_{k}$;

     \textbf{else}

     $\hspace{0.5cm}$ $\textup{\textbf{for}}$ ${t}$  $\textup{from}$ $1$ $\textup{to}$ $m_{k}$ $\textup{\textbf{do}}$

     $\hspace{1cm}$     $\textup{Find}$ $\textup{the}$ $\textup{least}$
                     $t_{k}$ $\textup{such}$ $\textup{that} $
                     $\nu_{t_{k}}\neq 0$;

     $\hspace{0.5cm}$ \textbf{end} \textbf{do};

     $\hspace{0.5cm}$ $\textup{\textbf{for}}$  ${t}$  $\textup{from}$     $t_k+1$  $\textup{to}$  $m_{k}$  $\textup{\textbf{do}}$

     $\hspace{2cm}$ $(  a_t,b_t) :=     (  a_t,b_t)-\cfrac{\nu_{t}}{\nu_{t_k}}(  a_{t_k},b_{t_k})$;

     $\hspace{0.5cm}$ $\textup{\textbf{end}}$ $\textup{\textbf{do}}$;

     $\hspace{1cm}$ $\mathcal{G}_{k+1}:=\Big\{(  a_1,b_1),\ldots,(
     a_{t_k-1},b_{t_k-1})$,$(  a_{t_k},b_{t_k})\cdot(x_1-y_{l,1}),\ldots,$

     $\hspace{3.5cm}$ $(  a_{t_k},b_{t_k})\cdot(x_n-y_{l,n})$,$(  a_{t_k+1},b_{t_k+1})$, $\ldots$,$(
a_{m_{k}},b_{m_{k}})\Big\}$

     $\hspace{1cm}$ $\mathcal{G}_{k+1}:= \textup{minimal~Gr\"{o}bner~basis} (\mathcal{G}_{k+1}) $

     \textbf{end} \textbf{if};

     $\textup{\textbf{return}}$ $\mathcal{G}_{k+1}$;

We must point out that here we do not require $\textbf{LC}((a,b))=1$
in minimal Gr\"{o}bner basis.

\begin{exmp}\label{eg3.1} Given the interpolation problem
\begin{table}[!htbp]
\centering
\begin{tabular}{|c|c|c|c|c|}
  \hline
  point & $f_i$ & $\frac{\partial}{\partial x} f$ & $\frac{\partial}{\partial y} f$ & $\frac{\partial^2}{\partial x y}f$ \\
  \hline
  (-1,2) & 2 &   &   &    \\
  \hline
  (1,1)  & 3 &   &   &   \\
  \hline
  (2,1)  & 4 & 5 & 2 &   \\
  \hline
  (3,2)  & 3 & 4 & 3 & 6 \\
  \hline
\end{tabular}
\caption{interpolation}
\end{table}

Fix the order $\prec_0$, using the Fitzpatrick algorithm\_RI, we can
compute the minimal Gr\"{o}bner basis of the submodule $M$:

$(a_1,b_1)=\big(\frac{1103}{14528}x^2-\frac{1367}{14528}xy-\frac{301}{7264}y^2+\frac{6713}{14528}x-\frac{959}{7264}y-1,~
-\frac{61}{908}y^2+\frac{3047}{14528}x+\frac{731}{14528}y-\frac{6335}{14528}\big),$

$(a_2,b_2)=\big(-\frac{19899}{314176}x^2+\frac{43619}{314176}xy-\frac{1999}{157088}y^2-\frac{153069}{314176}x+\frac{14059}{157088}y+1,~
\frac{122}{4909}xy+\frac{793}{19636}y^2-\frac{67507}{314176}x-\frac{19127}{314176}y+\frac{135787}{314176}\big),$

$(a_3,b_3)=\big(\frac{6973}{371696}x^2+\frac{61515}{371696}xy-\frac{16223}{185848}y^2-\frac{18057}{28592}x+\frac{4115}{185848}y+1,~
\frac{488}{23231}x^2+\frac{61}{1787}xy-\frac{89291}{371696}x-\frac{12399}{371696}y+\frac{141603}{371696}\big),$

$(a_4,b_4)=\big(-\frac{305}{12438}y^3-\frac{1519}{24876}x^2+\frac{2959}{24876}xy+\frac{1481}{12438}y^2-\frac{10769}{24876}x-\frac{673}{6219}y+1,~
\frac{122}{6219}xy+\frac{61}{1382}y^2-\frac{4847}{24876}x-\frac{1697}{24876}y+\frac{10027}{24876}\big),$

$(a_5,b_5)=\big(\frac{49}{988}xy^2-\frac{85}{494}y^3-\frac{6}{19}xy+\frac{214}{247}y^2+\frac{107}{247}x-y-\frac{22}{247},~
-\frac{15}{247}xy+\frac{61}{988}y^2+\frac{30}{247}x-\frac{12}{247}y-\frac{37}{247}\big),$

$(a_6,b_6)=\big(\frac{31}{474}x^2y-\frac{11}{158}xy^2+\frac{37}{237}y^3-\frac{31}{237}x^2+\frac{77}{158}xy-y^2-\frac{55}{79}x+\frac{78}{79}y+\frac{184}{237},~
\frac{1}{6}xy-\frac{55}{474}y^2-\frac{1}{3}x-\frac{5}{474}y+\frac{115}{237}\big),$

$(a_7,b_7)=\big(\frac{31}{978}x^3-\frac{11}{326}x^2y+\frac{37}{489}xy^2+\frac{12}{163}x^2-\frac{227}{978}xy-\frac{37}{163}y^2-\frac{679}{978}x+y+\frac{92}{163},~
\frac{79}{978}x^2-\frac{55}{978}xy-\frac{176}{489}x+\frac{55}{326}y+\frac{115}{326}\big)\Big\}$

Any weak interpolation $(a,b)$ have the form $$(  a,b)=c_1(
a_1,b_1)+\cdots+c_7(  a_7,b_7),$$ where $c_j\in \mathcal {P}$ \
$(j=1,\ldots,7)$ are free parameters. Choose $c_j$ properly such
that $b(Y_i)\neq 0,\ i=1,\ldots,L$, then we can get the
interpolation function
$$\cfrac{  a(X)}{b(X)}=\cfrac{c_1  a_1+\cdots +c_7  a_7}{c_1b_1+\cdots+c_7b_7}.$$
\end{exmp}

\section{Neville-like interpolation}\label{Neville}

Neville's algorithm is used for polynomial interpolation which was
derived by Eric Harold Neville. The algorithm aims at determining
the value of the interpolating polynomial at a single point $x$.
\cite{Key37} also derived a Neville type algorithm for univariable
rational interpolation.

In this section,  we present a Neville-like algorithm for
multivariate Cauchy interpolation based on Fitzpatrick
algorithm\_RI.

Given a set of $L$ distinct points $\{Y_1,\ldots,Y_L \}$, $Y_j\in
\mathbb{F}^n, j=1,\ldots,L$, and the corresponding values
$\{f_1,\ldots,f_L\}$, $f_j \in \mathbb{F}$, $j=1,\ldots,L$, we want
to determine the interpolating value at the point $Y_0$.

In this case, we know that $h_j:=f_j$, and $M_k:=\{( a,b)\in
M_{k-1}: b h_k - a=  0 \mod I_k\}=\{(  a,b)\in M_{k-1}: (b h_k-
a)\big|_{Y_k}= 0\}$, $k=1,\ldots,L$. In {Fitzpatrick algorithm\_RI},
if we can get the values $  W(i,j)=(b_i h_j- a_i)\big|_{Y_j}=
\nu_{i,j}$ and $(  a_i\big|_{(x_0,y_0)}, b_i\big|_{(x_0,y_0)})$
recursively without computing the weak interpolation $(a,b)$, then
we can determine the interpolating value at the point $Y_0$. It
means that using the present values, we can calculate the new values
of $W(i,j)$ and $( a_i\big|_{(x_0,y_0)}, b_i\big|_{(x_0,y_0)})$
without computing the weak interpolation $(a,b)$ when a new point is
added.  Based on this idea we get a Neville-Like algorithm for
Cauchy interpolation. For simplicity we will restrict ourselves to
the case $n=2$. Three and higher dimensional cases can be treated
similarly.

Fix the order $\prec_\xi$, and $y\prec x$. Let $\{(  a_1, b_1),
\ldots, ( a_{m_k}, b_{m_k})\}$ be the minimal Gr\"{o}bner basis of
$M_k$. Define $  W(i,j)=(b_i  h_j- a_i)\big|_{Y_j}$, $i=1,\ldots
m_k$, $j=1,\ldots,L$, $\vec  W(i,L+1)=(  a_i\big|_{(x_0,y_0)},
b_i\big|_{(x_0,y_0)})$, $\vec  {W}(i,L+2)=\textbf{LT}\big((
a_i,b_i)\big)$.

Define $\textbf{W}_i=\big(  W(i,1), \ldots, W(i,L), \vec W(i, L+1),
\vec W(i,L+2)\big)$.

We know that $\langle (  1,0)$, $( 0,1)\rangle$ is the Gr\"{o}bner
basis of $M_0=\mathcal{P}^{2}$.

First, using the Gr\"{o}bner basis of $M_0$, we compute

$\qquad$ $  W(1,j)=0\cdot  h_j-  1=-1,\ j=1,\ldots,L;$

$\qquad$ $  W(2,j)=1\cdot   h_j-  0=  f_j,\ j=1,\ldots,L.$

Let

$\qquad$ $\vec  W(1,L+1)=(  1,0),\ \vec  W(1,L+2)=\textbf{LT}((
1,0))=( 1,0),$

$\qquad$ $ \vec W(2,L+1)=(  0,1),\  \vec W(2,L+2)=\textbf{LT}((
0,1))=( 0,1)$.

We will compute the $\textbf{W}^{(k)}_{i}$ recursively by
$\textbf{W}_{1}$, $\textbf{W}_{2}$.

Define the matrix
$$\mathcal{W}_0:=\left(\begin{array}{c}
                                                           \textbf{W}_1 \\
                                                           \textbf{W}_{2}
                                                         \end{array}
\right)\overset{\Delta}{=}( \textbf{W}_1, \textbf{W}_{2})^T.$$

Let $\mathcal{W}_{k-1}=(\textbf{W}_1,\ldots,\textbf{W}_{m_{k}})^T$.
Denote by $\sharp \mathcal{W}_{k-1}$ the number of rows.

~\\ \textbf{Neville-like algorithm}:

     \textbf{Input}: $\mathcal{W}_0=(\textbf{W}_1,
     \textbf{W}_{2})^T$; $L$;

     \textbf{Output}:  $\mathcal{W}_{L}$;

\textbf{for} $k$ \textbf{from} $1$ \textbf{to} $L$ \textbf{do}

%
%

     $\hspace{1cm}$  $m_{k}:= \sharp \mathcal{W}_{k-1}$;

     $\hspace{1cm}$  $\textup{rearrange}$ $\textup{the}$ $\textup{elements}$ $\textup{of}$ $\mathcal{W}_{k-1}$ $\textup{so}$
     $\textup{that}$

     $\hspace{2.5cm}$               $\mathcal{W}_{k-1}=(\textbf{W}_1,\ldots,\textbf{W}_{m_{k}})^T$

     $\hspace{2.5cm}$               $\textup{and}$ $  {W}(1,L+2)\prec_{\xi}\cdots\prec_{\xi}   {W}(m_{k},L+2)$;




     $\hspace{1cm}$  \textbf{for} $i$ \textbf{from} $1$ \textbf{to}
     $m_{k}$ \textbf{do}

     $\hspace{1.5cm}$  $\textup{Find}$ $\textup{the}$ $\textup{least}$
                     $i_{0}$ $\textup{such}$ $\textup{that} $
                     $W(i_0,k)\neq 0$;

%
%
%
%
     $\hspace{1cm}$  \textbf{end} \textbf{do};
%
%


     $\hspace{1cm}$  \textbf{if}  $W(i,k)=0$ for all $i$  \textbf{then}

     $\hspace{2cm}$  $\mathcal{W}_{k}:=\mathcal{W}_{k-1}$

     $\hspace{1cm}$  \textbf{else}


     $\hspace{1.5cm}$ \textbf{for}  ${i}$  \textbf{from }    $i_0+1$  \textbf{to}  $m_{k}$ \textbf{do}

     $\hspace{2cm}$ \textbf{for} $j$ \textbf{from} 1 \textbf{to} $L$ \textbf{do}

     $\hspace{2.5cm}$
     $  W(i,j):=  W(i,j)-\cfrac{W(i,j_0)}{W(i_0,j_0)}\   W(i_0,j)$;

     $\hspace{2cm}$ \textbf{end} \textbf{do};

     $\hspace{2.cm}$
     $ \vec W(i,L+1):= \vec W(i,L+1)-\cfrac{W(i,j_0)}{W(i_0,j_0)}\  \vec W(i_0,L+1)$;

     $\hspace{2cm}$
     $ \vec W(i,L+2):= \vec W(i,L+2)$;

     $\hspace{1.5cm}$ \textbf{end} \textbf{do};

     $\hspace{2cm}$ \textbf{for} $j$ \textbf{from} 1 \textbf{to} $L$ \textbf{do}

     $\hspace{2.5cm}$
     $  W(m_{k}+1,j):=  W(i_0,j)\cdot (x_j-x_{j_0})$;

     $\hspace{2.5cm}$
     $  W(m_{k}+2,j):=  W(i_0,j)\cdot (y_j-y_{j_0})$;

     $\hspace{2cm}$ \textbf{end} \textbf{do};

     $\hspace{2.5cm}$
     $ \vec W(m_{k}+1,L+1):= \vec W(i_0,L+1)\cdot (x_0-x_{j_0})$;

     $\hspace{2.5cm}$
     $ \vec W(m_{k}+1,L+2):= \vec W(i_0,L+2)\cdot x$;

     $\hspace{2.5cm}$
     $ \vec W(m_{k}+2,L+1):= \vec W(i_0,L+1)\cdot (y_0-y_{j_0})$;

     $\hspace{2.5cm}$
     $ \vec W(m_{k}+2,L+2):= \vec W(i_0,L+2)\cdot y$;

     $\hspace{1.5cm}$ $\mathcal{W}_{k}:=\big(\textbf{W}_1,\ldots,
     \textbf{W}_{i_0-1}$,$\textbf{W}_{i_0+1},\ldots,\textbf{W}_{m_{k}+2}\big)^T$;

     $\hspace{1.5cm}$ $\mathcal{W}_k:=\textup{\textbf{Minimal}~\textbf{Gr\"{o}bner}~\textbf{basis}}(\mathcal{W}_k)$

     $\hspace{1cm}$ \textbf{end} \textbf{if};

     \textbf{end} \textbf{do};

     \textbf{return} $\mathcal{W}_{L}$;

~\\
$\textup{\textbf{Minimal}~\textbf{Gr\"{o}bner}~\textbf{basis}}(\mathcal{W})$

     \textbf{Input}: $\mathcal{W}$;

     \textbf{Output}:
     $\widetilde{\mathcal{W}}=(\textbf{W}_1,\ldots,\textbf{W}_{m})^T$
     with no $\vec W(i,L+2)$ is divisible by $\vec W(j,L+2)$ for
     $i\neq j$.

\begin{exmp}\label{eg4.1}

Given objective function $\ln(x^2+y^2)$, we will use the values at
the points $(1.75,1.75)$, $(2.25,1.75)$, $(1.75,2.25)$,
$(2.25,2.25)$, $(1.85,1.85)$, $(2.15,1.85)$, $(1.85,2.15)$, $(2.15,
2.15)$ to estimate the value at the point $(2,2)$.

Fix the order $\prec_0$, $L=8$, Neville-like algorithm outputs
$\mathcal{W}_8$:

\begin{table}[!htbp]
\centering
\begin{tabular}{|c|c|c|c|c|c|c|c|c|c|}
  \hline
  0.& 0.& 0.& 0.& 0.& 0.& 0.& 0.&($8.200\times10^{-7}$, $4.100\times10^{-7}$)& $(x^2, 0)$ \\
  0.& 0.& 0.& 0.& 0.& 0.& 0.& 0.&(12.45851102, 5.991319801)  &$(0, y^2)$ \\
  0.& 0.& 0.& 0.& 0.& 0.& 0.& 0.&(-37.67095363, -18.11585643)  & $(0, xy)$ \\
  0.& 0.& 0.& 0.& 0.& 0.& 0.& 0.&(37.67095324, 18.11585624)& $(0, x^2)$ \\
  0.& 0.& 0.& 0.& 0.& 0.& 0.& 0.&(-1.225410500, -0.5894116890)&$(y^3, 0)$ \\
  0.& 0.& 0.& 0.& 0.& 0.& 0.& 0.&(3.470473588, 1.668910238)& $(xy^2, 0)$ \\
  \hline
\end{tabular}
\caption{$\mathcal{W}_8$ for $\ln(x^2+y^2)$} \label{tb2}
\end{table}


From the $(8+1)$-th column of $\mathcal{W}_8$, we can see that each
of the vectors $\vec W(i,8+1)=( a_i\big|_{(2,2)},b_i\big|_{(2,2)})$,
$i=1,\ldots,6$, gives an approximate value
$\cfrac{a_i|_{(2,2)}}{b_i|_{(2,2)}}$ of $\ln(2^2+2^2)$.

Here we choose
$$\cfrac{\sum_{i=1}^{m_k} \textbf{sgn}(b_i\big|_{Y_0})\cdot
 a_i\big|_{Y_0}}{\sum_{i=1}^{m_k}
\textbf{sgn}(b_i\big|_{Y_0})\cdot b_i\big|_{Y_0}},$$ as our
estimation value(see Table \ref{tb3}), where $Y_0=(2,2)$,
$\textbf{sgn}(x)$ satisfies: if $x\geq 0$, $\textbf{sgn}(x)=1$, else
$\textbf{sgn}(x)=-1$.
\begin{table}[!htbp]
\centering
\begin{tabular}{|c|c|c|c|}
  \hline
  $i$& $(x_i,y_i)$  &  $\ln(x_i^2+y_i^2)$ & interpolating value of $\mathcal {W}_i$ \\
  1 & (1.75,1.75)   & 1.812378756 & 1.312378756 \\
  2 & (2.25,1.75)   & 2.094945728 & 1.812378756 \\
  3 & (1.75,2.25)   & 2.094945728 & 2.122484930  \\
  4 & (2.25,2.25)   & 2.315007613 & 2.107686660  \\
  5 & (1.85,1.85)  & 1.923518459 & 2.082067864 \\
  6 & (2.15,1.85)  & 2.085050780 & 2.082067864 \\
  7 & (1.85,2.15)  & 2.085050780 & 2.079431546\\
  8 & (2.15,2.15)  & 2.224082865 & 2.079439873\\
  \hline
\end{tabular}
\caption{estimation value of $\ln(2^2+2^2)$}\label{tb3}
\end{table}
Actually $\ln(2^2+2^2)=2.079441542$.
\end{exmp}

\begin{exmp}\label{eg4.3}

Given the values of $\sqrt{1-x^2-y^2}$ at the points $(0.45,0.5)$,
$(0.55,0.45)$, $(0.45,0.55)$, $(0.55,0.55)$, $(0.5,0.45)$,
$(0.5,0.55)$, $(0.45,0.55)$, $(0.55$, $0.5)$. Fix the order
$\prec_0$, we estimate the value of $\sqrt{1-(0.5)^2-(0.5)^2}$ by
$\frac{\sum_{i=1}^{m_k} \textbf{sgn}(b_i\big|_{(0.5,0.5)})\cdot
 a_i\big|_{(0.5,0.5)}}{\sum_{i=1}^{m_k}
\textbf{sgn}(b_i\big|_{(0.5,0.5)})\cdot b_i\big|_{(0.5,0.5)}}$ (see
Table \ref{tb4}).

\begin{table}
\centering
\begin{tabular}{|c|c|c|c|}
  \hline
  $i$ & $Y_i$&  $\sqrt{1-(x_i)^2-(y_i)^2}$& interpolating value of $\mathcal {W}_i$ \\
  1 & (0.45,0.45) & 0.7713624310 & 0.6713624310  \\
  2 & (0.55,0.45) & 0.7035623640 & 0.7078673362 \\
  3 & (0.45,0.55) & 0.7035623640 & 0.7035623636 \\
  4 & (0.55,0.55) & 0.6284902545 & 0.7035623639 \\
  5 & (0.5,0.45)  & 0.7399324293 & 0.7035623639 \\
  6 & (0.5,0.55)  & 0.6689544080 & 0.7047928585 \\
  7 & (0.45,0.5)  & 0.7399324293 & 0.7071486038 \\
  8 & (0.55,0.5)  & 0.6689544080 & 0.7071187945 \\
  \hline
\end{tabular}
\caption{estimation value of $\sqrt{1-(0.5)^2-(0.5)^2}$}\label{tb4}
\end{table}
Actually $\sqrt{1-(0.5)^2-(0.5)^2}=0.7071067812$.
\end{exmp}

\begin{exmp}\label{eg4.4}

Given the values of $\exp(x^2+y)$ at the points $(2, 2.95)$, $(2,
3.05)$, $(1.95, 3)$, $(2.05, 3)$, $(1.975, 2.975)$, $(1.975,
3.025)$, $(2.025, 2.975)$, $(2.025, 3.025)$. Fix the order
$\prec_0$, we still choose $\frac{\sum_{i=1}^{m_k}
\textbf{sgn}(b_i\big|_{(2,3)})\cdot
 a_i\big|_{(2,3)}}{\sum_{i=1}^{m_k}
\textbf{sgn}(b_i\big|_{(2,3)})\cdot b_i\big|_{(2,3)}}$ as our
estimation value(see Table \ref{tb5}).

\begin{table}
\centering
\begin{tabular}{|c|c|c|c|}
  \hline
  $i$ & $Y_i$&  $\exp(x_i^2+y_i)$ & interpolating value of $\mathcal {W}_i$ \\
  1 & (2, 2.95) & 1043.149728 & 1043.099728  \\
  2 & (2, 3.05) & 1152.858743 & 1044.131824 \\
  3 & (1.95, 3) & 900.0947180 & 1043.425504 \\
  4 & (2.05, 3) & 1342.783531 & 1102.658424 \\
  5 & (1.975, 2.975)  & 968.3804142 & 1097.459656 \\
  6 & (1.975, 3.025)  & 1018.030340 & 1096.945601 \\
  7 & (2.025, 2.975)  & 1182.782509 & 1096.552830 \\
  8 & (2.025, 3.025)  & 1243.425065 & 1096.660126 \\
  \hline
\end{tabular}
\caption{estimation value of $\exp(2^2+3)$}\label{tb5}
\end{table}
Actually $\exp(2^2+3)=1096.633158$.
\end{exmp}

\section{Conclusion}\label{Conclusion}

In this paper, we apply the Fitzpatrick algorithm to osculatory
rational interpolation, and get the parametric solution of all the
interpolation functions with the given complexity .

For Cauchy interpolation, we present a Neville-like algorithm to
determine the value of interpolating function at a single point
without computing the rational interpolation function (several
points can be treated  similarly). Since each of the vectors $\vec
W(i,L+1)=(a_i\big|_{Y_0},b_i\big|_{Y_0})$ gives an approximate
value, we choose $$\cfrac{\sum_{i=1}^{m_k}
\textbf{sgn}(b_i\big|_{Y_0})\cdot
 a_i\big|_{Y_0}}{\sum_{i=1}^{m_k}
\textbf{sgn}(b_i\big|_{Y_0})\cdot b_i\big|_{Y_0}},$$ as our
estimation value. From the examples we can see that the Neville-like
algorithm is effective.

\bibliographystyle{elsarticle-num}
\bibliography{ref}

\end{document}